\documentclass[11pt]{article}

\usepackage{a4wide}
\usepackage{enumerate}
\usepackage[greek, frenchb, english]{babel}
\usepackage[T1]{fontenc}        
\usepackage{mathrsfs} 
\usepackage{amssymb}
\usepackage{amsmath}
\usepackage{amsfonts}
\usepackage[dvips]{graphicx}
\usepackage{color} %\textcolor{orange{texte \`a colorier}
\usepackage[colorlinks=true, urlcolor=blue,linkcolor=blue, citecolor=blue]{hyperref}
\usepackage{aeguill}
\usepackage{rotating}
\usepackage{units}
\usepackage{graphicx}

\allowdisplaybreaks
\date{}

\newtheorem{theo}{Theorem}[section]
\newtheorem{df}[theo]{Definition}%[section]
\newtheorem{lemme}[theo]{Lemma}%[section]
\newtheorem{Rem}[theo]{Remark}%[section]

\newtheorem{prop}[theo]{Proposition}

\newcommand{\Proof}{{\bf Proof:}}
\newcommand{\CQFD}{\hfill $\square$}
\newcommand{\ind}{\mathbf{1}}

\numberwithin{equation}{section}

\def\real{\Bbb{R}}
\def\rdr+{\Bbb{R}^d\times\Bbb{R}_+}
\def\rrr{\Bbb{R}^d\times\Bbb{R}_+\times\Bbb{R}}
\def\rd{\Bbb{R}^d}
\def\r+d{\Bbb{R}_+}

\def\ee{\Bbb{E}}

\def\disp{\displaystyle}

\title{Macroscopic analysis of shot-noise Cox random balls}

\author{
Adrien Clarenne\footnotemark[1] \
}

\begin{document}

\maketitle

\footnotetext[1]{Univ Rennes, CNRS, IRMAR - UMR 6625, F-35000 Rennes, France. Email: adrien.clarenne@univ-rennes1.fr}

\begin{abstract}
In this paper, we consider a cluster model of weighted Euclidean random balls generated by a shot-noise Cox process. It is an example of cluster point process. We perform a scaling on the model by shrinking the radii of the balls and compensate this effect by increasing the (mean) number of balls in each cluster, or/and increasing the (mean) number of clusters. We consider two different scenarios, say a local and a global scenarios. Heuristically, in the first scenario, we focus on the mean number of large balls in a cluster while in the second one, we focus on the global mean number of large balls in the model. According to the different scenarios, the cluster structure can persist at the limit or disappear.\\\\
{\bf Keywords:}
Cox processes, 
random balls,
random fields, 
limit theorem.
\\
{\bf MSC Classification:} 
60G55, %point process
60F05, %limit theorem
60G60. %random fields
\end{abstract}

%%%%%%%%%%%%%%%%%%%%%%%%%%%%%%%%%%%%%%%%%%%%%%%%%%%%%%%%%%%%%%%%%%%%%%%%%%%%%%%%%%%%%%%%%%%%%%%%%%%%

\section*{Introduction}

We consider a model of weighted Euclidean random balls in $\rd$ generated by a shot-noise Cox process as follows. The centers of the balls are generated by a shot-noise Cox process $Z$ and this point process is marked twice, first by random variables $r$ with density probability $f$, seen as the radii of the balls, and second by a mark seen as the weights of the balls generated by a probability measure $G$. The marks are all independent and independent also of $Z$. The resulting marked point process is still a Cox process $C$ but on $\rd\times\real_+\times\real$. \\
It is an example of cluster Poisson process where the centers of the clusters are drawn by a Poisson point process $\Phi$ on $\rd$. Since a Cox process can be seen as a Poisson point process with a random intensity, the model under study is actually a randomized version of the Poissonian random balls model. When the intensity is deterministic, we recover the usual Poissonian model which has been studied in several papers since its introduction in \cite{KLNS2007}, see also \cite{BEK2010, BD2009} for generalizations of \cite{KLNS2007} with weights and/or zoom-in and zoom-out. See also \cite{BCG2018} for a determinantal random ball model beyond the Poissonian setting.\\
The class of Cox point process is one of the most used among the cluster models, because it can represent random constraint of a field, for example a random heterogeneity (see \cite{Moller2003, MT2005} for more details on shot-noise Cox processes). In dimension 1, a random balls model can be interpreted as a model for the study of communication network (see \cite{MRRS2002}) or power consumption for example. In this situation, the centers $x$ of the balls are interpreted as the date of the connection and the radius $r$ the duration of the connection. In the Poissonian case, the dates of connections are placed uniformly in time. In the determinantal case, the connection dates are not too close to each other. Here, in the shot-noise Cox model, we represent the situation where there are peaks of connections, for example in the morning or in the evening when people are at home. We are typically in a cluster situation. In dimension~2, we can interpret the model as a wireless network, where $x$ is the location of a transmitter and $r$ is its range of transmission. In the Poissonian case, the transmitters are uniformly distributed in the space. In the shot-noise Cox model, the cluster structure implies that the antennas are highly concentrated in some places and sparse at others which is indeed the case in some real situation (think about a city with no antennas on lakes, river or some special infrastructure like schools). \\
In the following, our macroscopic analysis is driven as follow. We perform a scaling in this model by first shrinking the radii; to compensate this effect, we rescale the shot-noise Cox process $Z$ that generated the centers of the balls. In contrast with the Poissonian (see \cite{BEK2010}) or the determinantal case (see \cite{BCG2018}), where there is just one level of randomness, the global location of the balls, here we have two levels of randomness with the collection of cluster and within each cluster. As a consequence, this additional level of randomness makes it possible to study many more different asymptotics behaviours in this model. In practice, we dispose of two levels of scaling with the (mean) number of clusters and the (mean) number of balls in each cluster.\\
In the Poissonian or determinantal cases, a key quantity appearing to drives the fluctuations is the mean number of large balls. In the setting of Cox process, two different scenarios, say a local and a global one, are possible because two distinct quantities can drive the fluctuations.
The first one is the mean number of large balls in each cluster, and we will refer to it as the local scaling. In this scenario, we do not rescale the mean number of clusters. Roughly speaking, each cluster is a Poissonian ball model whose asymptotics are well known from \cite{BEK2010, BD2009} and the whole limit of the Cox model is then a mixture of the limit random fields obtained, and so we obtain a randomized version of the Poissonian results from \cite{BEK2010, BD2009}.\\
The second scenario will be refer to as the global scaling: in contrast to the first one where we only focus on the mean number of balls in each cluster, in this scenario we focus on the global mean number of large balls in the model. This situation is analogous to the Poissonian case, and we recover the same three regimes as therein, with a disappearance of the cluster structure.\\
The document is organized as follow. In Section \ref{sec-model}, we give a detailed description of the model under study. In Section \ref{sect-scaled model}, we introduce the rescaled model and the object of interest in the paper. In Section \ref{sec:asympt}, we give the main results in two subsections dealing with the two scenarios described above. Finally, some general results about Cox processes are given in Appendix~\ref{appendix}. In particular, we refer to this Appendix for any reader unfamiliar with Cox processes, see Definition \ref{def-cox}.

\section{Shot-noise Cox random ball model}\label{sec-model}

\noindent We now describe mathematically the model under study. We consider a collection of Euclidean balls $B(x,r)=\left\lbrace y\in\rd : \|y-x\|\leq r \right\rbrace$ whose centers $x$ and radii $r$ are generated by a Cox process on $\rdr+$. To generate the balls, we first consider a shot-noise Cox process $D$ on $\rd$ directed by a random measure 
\begin{equation}\label{def:Z}
 Z(dx)=\sum_{y\in\Phi}k(x,y)dx
\end{equation}
that generates the centers of the balls. Here, $\Phi$ is a Poisson point process on $\rd$ with intensity the Lebesgue measure on $\rd$ and $k$ a positive function on $\rd\times\rd$ satisfying
\begin{equation}\label{k-bornee}
\|k\|_{\infty}=\underset{(x,y)\in(\real^d)^2}{\sup}k(x,y)<+\infty
\end{equation}
and for all $x\in\rd$, 
\begin{equation}\label{eq:k-intensity}
 \int_{\rd}k(x,y)dy=1.
\end{equation}
This is an example of a Poisson cluster process, where $\Phi$ is the base point process drawing the centers $c$ of the clusters $X_c$, and $k(c,\cdot)$ is the density intensity of the Poisson process $X_c$.\\
To each center $x$, we attach two marks $r$ (positive) and $m$, interpreted respectively as the radius of the ball and as the weight of the ball. These radii and weights are each identically distributed according respectively to a density $f$ on $\r+d$ and a probability distribution $G$ on $\real$. We obtain a Cox process $C$ on $\rrr$ directed by the random measure 
\begin{equation}\label{int-cox-poids}
\Lambda(dx,dr,dm)=Z(dx)f(r)drG(dm). 
\end{equation} 
Like in the previous studies of random ball model \cite{BEK2010, BD2009, BCG2018}, we focus on the following measure-indexed random field given for a measure $\mu$ on $\rd$ by
\begin{equation}\label{expr:M}
M(\mu)=\int_{\rrr}m\mu(B(x,r))C(dx,dr,dm).
\end{equation} 
To ensure that the quantity \eqref{expr:M} is indeed well defined, we actually restrain the study of $M$ to finite signed measures $\mu$, i.e. $\mu\in \mathcal{Z}(\rd):=\left\lbrace \mu : \vert\mu\vert(\rd)<+\infty \right\rbrace$ and we denote $\|\mu\|=\vert\mu\vert(\rd)$ for the total variation of $\mu$. We also assume that the distribution $G$ of the marks belongs to the normal domain of attraction of the $\alpha$-stable distribution $S_\alpha (\sigma,b,0)$ with $\alpha \in (1,2]$. Here, following the terminology of \cite{ST1994}, $\sigma$ is a scale parameter and $b$ is a skewness parameter while the translation parameter is zero. Since $\alpha >1$, we have
\begin{equation}\label{esp:G finie}
\int_{\real}\vert m\vert G(dm)<+\infty.
\end{equation}
Moreover, like in previous studies of random ball models \cite{BEK2010,BD2009,BCG2018}, we assume also the following power-law hypothesis on the radius behaviour: for $d<\beta <\alpha d$,
\begin{equation}\label{hyp:f}
f(r)\underset{r\rightarrow +\infty}{\sim}\frac{C_\beta}{r^{\beta +1}}, \hspace{0.2cm} f(r)\leq \frac{C_0}{r^{\beta +1}},
\end{equation}
for finite positive constant $C_0$ and $C_\beta$.\\
Observe in particular that since $\beta>d$, the mean volume of the ball is finite:
\begin{equation}\label{eq:vol fini}
\int_{\rd} r^df(r)dr<+\infty.\\
\end{equation}
For this Cox process, we show now that the quantity \eqref{expr:M} is indeed well defined for $\mu \in\mathcal{Z}(\rd)$: using Fubini-Tonelli theorem, we have

\begin{align*}
\ee\left[|M(\mu)| \right] &= \ee\left[\left|\int_{\rrr}m\mu(B(x,r))C(dx,dr,dm)\right|\right]\\
&=\ee\left[\ee\left[\left|\int_{\rrr}m\mu(B(x,r))C(dx,dr,dm)\right|\bigg|\Lambda\right]\right]\\
&=\ee\left[\left|\int_{\rrr}m\mu(B(x,r))\Lambda(dx,dr,dm)\right|\right]\\
&\leq\ee\left[\int_{\rd}\int_{\rrr}\vert m\vert \vert \mu(B(x,r))\vert k(x,y)dxf(r)drG(dm) \Phi(dy)\right]\\
&=\int_{\rd}\int_{\rrr}\vert m\vert \vert \mu(B(x,r))\vert k(x,y)dxf(r)drG(dm) dy\\
&=\left(\int_{\real}\vert m\vert G(dm)\right)\left(\int_{\rdr+}\vert \mu(B(x,r))\vert dxf(r)dr \right)\\
&\leq v_d \|\mu\| \int_{\real}\vert m\vert G(dm)\int_{\r+d}r^df(r)dr <+\infty
\end{align*}
thanks to \eqref{esp:G finie} and \eqref{eq:vol fini}, where $v_d$ stands for the volume of the unit Euclidean ball of $\rd$.\\

%%%%%%%%%%%%%%%%%%%%%%%%%%%%%%%%%%%%%%%%%%%%%%%%%%%%%%%%%%%%%%%%%%%%%%%%%%%%%%%%%%%%%%%%%%%%%%%%%%%%%
%%%%%%%%%%%%%%%%%%%%%%%%%%%%%%%%%%%%%%%%%%%%%%%%%%%%%%%%%%%%%%%%%%%%%%%%%%%%%%%%%%%%%%%%%%%%%%%%%%%%%

\section{Rescaled model}\label{sect-scaled model}

Our asymptotics are obtained by zooming-out in the model. This is obtained by performing a scaling in the model. To that purpose, we introduce a zooming-out rate $\rho \in ]0,1[$ to shrink the radii of the balls. In order to compensate the shrinking of the balls, the random intensity measure $Z(dx)$ in \eqref{def:Z} is changed into $Z_\rho(dx)$ with
\begin{equation*}
 Z_\rho(dx)=\sum_{y\in\Phi_\rho}k_\rho(x,y)dx,
\end{equation*} 
where $\Phi_\rho$ is a Poisson point process on $\rd$ with intensity measure $\kappa(\rho)dy$ and $k_\rho$ is a density kernel satisfying:
\begin{equation}\label{hyp-krho}
k_\rho(x,y)\underset{\rho\rightarrow 0}{\sim}\lambda(\rho)k(x,y)~\text{and}~k_\rho(x,y)\leq C_1\lambda(\rho)k(x,y).
\end{equation}
The parameters $\kappa(\rho)$ and $\lambda(\rho)$ will compensate the zooming-out effect (i.e. the balls become smaller when $\rho$ goes to 0). Heuristically, $\lambda(\rho)$ gives the order of the mean number of balls in a cluster and $\kappa(\rho)$ gives that of the clusters.
\\
Let $C_\rho$ be the resulting scaled version of the Cox process $C$ on $\rrr$ given in \eqref{int-cox-poids} i.e. $C_\rho$ is directed by the scaled intensity measure 
\begin{equation*}\label{int-cox-scaled}
\Lambda_\rho(dx,dr,dm)=Z_\rho(dx)f\left(\frac{r}{\rho}\right)\frac{dr}{\rho }G(dm)
\end{equation*}
(it is an analogous scaling as in \cite{BEK2010,BCG2018,Gobard_phd}) and we are interested in the corresponding rescaled quantity as in \eqref{expr:M}:
\begin{equation}\label{objet-rescaled}
M_\rho(\mu)=\int_{\rrr}m\mu(B(x,r))C_\rho(dx,dr,dm).
\end{equation}
In the Poissonian case, the fluctuations of $M_\rho(\mu)$ with respect to its mean value are investigated. Since Cox processes can be seen as Poisson processes but with random intensity measure, it is relevant in the present Cox setting to study the fluctuations of $M_\rho(\mu)$ in \eqref{objet-rescaled} with respect to its conditional mean. In this situation, the centering is thus not deterministic and we investigate the limit in law when $\rho\rightarrow 0$ of
\begin{equation}\label{expr:fluct}
\frac{M_\rho(\mu)-\ee\left[M_\rho(\mu)\left.\right|\Lambda_\rho\right]}{n(\rho)}
\end{equation}
for a proper normalization $n(\rho)$. \\

The fluctuations of $M_\rho(\mu)$ given in \eqref{expr:fluct} are ruled by the mean number of large balls in the model. By large balls, we mean balls, say, with radius larger than one and containing the origin 0. Let compute this key quantity. Setting $\# A$ for the cardinal of $A$,

\begin{align*}
&\ee\bigg[\#\left\lbrace (x,r,m)\in C_\rho \big| 0\in B(x,r),~r>1\right\rbrace \bigg]\\
&\hspace*{4cm}=\ee \bigg[\int_{\rrr}\ind_{\left\lbrace 0\in B(x,r),~r>1\right\rbrace}C_\rho(dx,dr,dm)   \bigg]\\
&\hspace*{4cm}=\ee \bigg[\ee \bigg[\int_{\rrr}\ind_{\left\lbrace 0\in B(x,r),~r>1\right\rbrace}C_\rho(dx,dr,dm)\bigg| \Lambda_\rho \bigg] \bigg]\\
&\hspace*{4cm}=\ee \bigg[\int_{\rrr}\ind_{\left\lbrace 0\in B(x,r),~r>1\right\rbrace}\Lambda_\rho(dx,dr,dm)\bigg]\\
&\hspace*{4cm}=\ee \bigg[\int_{\rrr}\ind_{\left\lbrace 0\in B(x,r),~r>1\right\rbrace}Z_\rho(dx)f\left(\frac{r}{\rho}\right)\frac{dr}{\rho }G(dm)\bigg]\\
&\hspace*{4cm}=\ee \bigg[\int_{\rd} \int_{\rrr}\ind_{\left\lbrace 0\in B(x,r),~r>1\right\rbrace}k_\rho(x,y)dx f\left(\frac{r}{\rho}\right)\frac{dr}{\rho }G(dm) \Phi_\rho(dy)\bigg]\\
&\hspace*{4cm}=\int_{\rd} \int_{\rrr}\ind_{\left\lbrace 0\in B(x,r),~r>1\right\rbrace}k_\rho(x,y)dx f\left(\frac{r}{\rho}\right)\frac{dr}{\rho }G(dm) \kappa(\rho)dy\\
&\hspace*{4cm}=\kappa(\rho)\int_{\rd} \int_{\rd}\int_1^{+\infty}\ind_{\left\lbrace x\in B(0,r)\right\rbrace}k_\rho(x,y) f\left(\frac{r}{\rho}\right) dxdy\frac{dr}{\rho }.
\end{align*}
The dominated convergence theorem gives the behaviour of this last integral:\\
\begin{itemize}
\item $\disp \underset{\rho\rightarrow 0}{\lim} ~\frac{1}{\lambda (\rho)\rho^\beta}\times\ind_{x\in B(0,r)}k_\rho(x,y) f\left(\frac{r}{\rho}\right) \frac{1}{\rho }= \ind_{x\in B(0,r)}k(x,y)C_\beta r^{-\beta-1}$
\item $\disp\bigg|\frac{1}{\lambda (\rho)\rho^\beta}\times\ind_{x\in B(0,r)}k_\rho(x,y) f\left(\frac{r}{\rho}\right) \frac{1}{\rho }\bigg|\leq  \ind_{x\in B(0,r)}C_1 k(x,y)C_0 r^{-\beta-1}$
\end{itemize}
which is independent of the parameter $\rho$ and integrable on $\rd\times\rd\times ]1,+\infty[$ since:
\begin{align*}
\int_{\rd} \int_{\rd}\int_1^{+\infty}\ind_{x\in B(0,r)}C_1 k(x,y)C_0 r^{-\beta-1}dxdydr&=C_0C_1\int_{\rd} \int_1^{+\infty}\ind_{x\in B(0,r)} r^{-\beta-1}dxdr\\
&=C_0C_1 v_d \int_1^{+\infty} r^{d-\beta-1}dxdr \\
&=\frac{C_0C_1 v_d }{\beta -d} <+\infty.
\end{align*}
Thus we have:
\[\underset{\rho\rightarrow 0}{\lim}\int_{\rd} \int_{\rd}\int_1^{+\infty}\frac{1}{\lambda (\rho)\rho^\beta}\times\ind_{x\in B(0,r)}k_\rho(x,y) f\left(\frac{r}{\rho}\right) \frac{1}{\rho }dxdydr=\frac{C_\beta v_d }{\beta -d},\]
which finally gives:
\begin{equation}\label{nb-grosse-boule}
\ee\bigg[\#\left\lbrace (x,r,m)\in C_\rho \big| 0\in B(x,r),~r>1\right\rbrace \bigg]\underset{\rho\to 0}{\sim}\frac{C_\beta v_d }{\beta -d}\kappa(\rho)\lambda (\rho)\rho^\beta.
\end{equation}
This computation shows that the mean number of large balls in the model is of order $\kappa(\rho)\lambda (\rho)\rho^\beta$. We can also interpreted this result as follows: $\kappa(\rho)$ represents the mean number of clusters in the model, so the mean number of large balls in each cluster is of order $\lambda (\rho)\rho^\beta$. These two interpretations will give two different studies of the model, as explained at the beginning of Section \ref{sec:asympt}.\\

The limit in law of \eqref{expr:fluct} when $\rho$ goes to 0 will be identified by the asymptotics of its characteristic function for which we dispose of the expression \eqref{eq:laplace-centree} for the shot-noise Cox process given in the Appendix. Applying Proposition \ref{prop:laplace-cox-centree} and Fubini-Tonelli theorem: 
\begin{align*}
&\ee\left[\exp\left(i\theta\left(\frac{M_\rho(\mu)-\ee\left[M_\rho(\mu)\left.\right|\Lambda_\rho\right]}{n(\rho)} \right)\right) \right]\\
&\hspace*{1.2cm}=\ee \left[\exp\left(\int_{\rrr} \psi\left(\frac{\theta m\mu(B(x,r))}{n(\rho)}\right) \Lambda_\rho(dx,dr,dm)\right) \right]\\
&\hspace*{1.2cm}=\ee \left[\exp\left(\int_{\rdr+} \psi_G\left(\frac{\theta \mu(B(x,r))}{n(\rho)}\right) Z_\rho(x)dxf\left(\frac{r}{\rho}\right) \frac{dr}{\rho}\right) \right]\\
&\hspace*{1.2cm}=\ee \left[\exp\left(\int_{\rdr+}\int_{\rd} \psi_G\left(\frac{\theta \mu(B(x,r))}{n(\rho)}\right)k_\rho(x,y)\Phi_\rho(dy) dxf\left(\frac{r}{\rho}\right) \frac{dr}{\rho}\right) \right]\\
&\hspace*{1.2cm}=\ee \left[\exp\left(\int_{\rd}\left(\int_{\rdr+} \psi_G\left(\frac{\theta \mu(B(x,r))}{n(\rho)}\right) k_\rho(x,y)dxf\left(\frac{r}{\rho}\right) \frac{dr}{\rho}\right)\Phi_\rho(dy)\right) \right],
\end{align*}
where $\displaystyle \psi_G (u)=\int_\real \psi(mu)G(dm)$ and $\psi(x)=e^{ix}-1-ix$.\\
In order to apply Proposition \ref{prop:laplace-poisson}, we check that \eqref{condition-poisson} holds true, with:
 $$\displaystyle \widetilde{g}(y)=-i\int_{\rdr+} \psi_G\left(\frac{\theta  \mu(B(x,r))}{n(\rho)}\right) k_\rho(x,y)dxf\left(\frac{r}{\rho}\right) \frac{dr}{\rho}.$$
We have $\left|\psi(u)\right|\leq 2|u|$ (see Lemma 1 in \cite{KLNS2007} for details) which gives the following control for $\psi_G$:
\begin{equation}\label{ineg-psig}
\vert\psi_G (u)\vert\leq 2\left(\int_\real \vert m \vert G(dm)\right) \vert u\vert.
\end{equation}
Also using \eqref{hyp-krho}, we have the following domination:

\begin{align*}
&\int_{\rd}(1\wedge|\widetilde{g}(y)|)\kappa(\rho)dy \leq\kappa(\rho)\int_{\rd}|\widetilde{g}(y)|dy \\
&\hspace*{0.7cm}=\kappa(\rho)\int_{\rd}\left|\int_{\rdr+} \psi_G\left(\frac{\theta \mu(B(x,r))}{n(\rho)}\right) k_\rho(x,y)dxf\left(\frac{r}{\rho}\right) \frac{dr}{\rho}\right|dy \\
&\hspace*{0.7cm}\leq \kappa(\rho)\int_{\rd}\int_{\rdr+} \left|\psi_G\left(\frac{\theta \mu(B(x,r))}{n(\rho)}\right)\right| k_\rho(x,y)dxf\left(\frac{r}{\rho}\right) \frac{dr}{\rho}dy\\
&\hspace*{0.7cm}\leq 2\kappa(\rho)\left(\int_\real \vert m \vert G(dm)\right)\int_{\rd}\int_{\rdr+} \frac{\vert\theta \mu(B(x,r))\vert}{n(\rho)}k_\rho(x,y)dxf\left(\frac{r}{\rho}\right) \frac{dr}{\rho}dy \\
&\hspace*{0.7cm}\leq \frac{2 C_1 \vert\theta\vert  \lambda (\rho)\kappa(\rho)}{n(\rho)}\left(\int_\real \vert m \vert G(dm)\right) \int_{\rd}\int_{\rdr+}  \vert\mu(B(x,r))\vert k(x,y)dxf\left(\frac{r}{\rho}\right) \frac{dr}{\rho}dy.
\end{align*}
By Fubini-Tonelli theorem, using (\ref{eq:k-intensity}), we get:
\begin{align*}
&\int_{\rd}\int_{\rdr+}  \vert \mu(B(x,r))\vert k(x,y)dxf\left(\frac{r}{\rho}\right) \frac{dr}{\rho}dy\\
&\hspace*{5.5cm}=\int_{\rdr+} \vert \mu(B(x,r))\vert \overbrace{\left(\int_{\rd} k(x,y) dy \right)}^{=1}  dxf\left(\frac{r}{\rho}\right) \frac{dr}{\rho} \\
&\hspace*{5.5cm}=\int_{\rdr+}\vert  \mu(B(x,r))\vert dxf\left(\frac{r}{\rho}\right) \frac{dr}{\rho}\\
&\hspace*{5.5cm}\leq C_1 \|\mu\|v_d \rho^d \int_{\real_+} r^df(r)dr.
\end{align*}
Finally:
\begin{align*}
\int_{\rd}|g(y)|\kappa(\rho)dy&\leq\frac{2 C_1 \vert\theta\vert  \lambda (\rho)\kappa(\rho)}{n(\rho)}\left(\int_\real \vert m \vert G(dm)\right)\int_{\rdr+} \vert \mu(B(x,r))\vert dxf\left(\frac{r}{\rho}\right) \frac{dr}{\rho} \\
&=\frac{2C_1 \vert \theta\vert  \|\mu\|v_d \lambda (\rho)\rho^d\kappa(\rho)}{n(\rho)}\left(\int_{\real}\vert m\vert G(dm)\right)\int_{\real_+}r^d f(r)dr <+\infty.
\end{align*}
As a consequence Proposition \ref{prop:laplace-poisson} applies and gives the characteristic function of \eqref{expr:fluct}:

\begin{prop}
Assume \eqref{eq:k-intensity}, \eqref{hyp:f} and \eqref{hyp-krho} and let $\mu\in \mathcal{Z}(\rd)$. Then the characteristic function of \eqref{expr:fluct} is given by
\begin{align}\label{laplace-shot-noise}
&\ee\left[\exp\left(i\theta\left(\frac{M_\rho(\mu)-\ee\left[M_\rho(\mu)\left.\right|\Lambda_\rho\right]}{n(\rho)} \right)\right) \right]\nonumber\\
&\hspace{-0.1cm}=\exp\left(-\int_{\rd} \left(1-\exp\left(\int_{\rdr+} \psi_G\left(\frac{\theta \mu(B(x,r))}{n(\rho)}\right)k_\rho(x,y) dxf\left(\frac{r}{\rho}\right) \frac{dr}{\rho}\right) \right)\kappa(\rho) dy \right).
\end{align}
\end{prop}

In order to investigate the behaviour of \eqref{expr:fluct} when $\rho \to 0$, i.e. the limit of \eqref{laplace-shot-noise}, it is necessary to consider a restricted class of measures $\mu$ that we introduce now (see also \cite{BD2009, BCG2018}).
\begin{df}
\label{def:Mbeta}
The set ${\cal M}_{\alpha,\beta}$ consists of signed measures $\mu \in {\cal Z}(\real^d)$ such that there exists two real numbers $p$ and $q$ with $0<p<\beta< q\leq 2d$ and a positive constant $C_\mu$ such that 
\begin{equation}
\label{eq:Mbeta}
\int_{\real^d}\vert\mu\big(B(x,r)\big)\vert^\alpha\ dx\leq C_\mu\big( r^{p}\wedge r^{q}\big),
\end{equation}
where $a\wedge b=\min(a,b)$.
\end{df}
%vide quand q>2d: cf BEK, p.1116.
The control in \eqref{eq:Mbeta} by both $r^p$ and $r^q$ is required to ensure that the integral in Proposition~\ref{prop:Mcal}-(i) below is indeed well defined. This integral is actually of constant use in our argument, so the introduction of the space ${\cal M}_{\alpha,\beta}$ is crucial. This definition is reminiscent of ${\cal M}_{2,\beta}$ in \cite{BD2009}. In particular, absolutely continuous measures with respect to the Lebesgue measure, with density $\varphi\in L^1(\real^d)\cap L^\alpha(\real^d)$, do belong to ${\cal M}_{\alpha,\beta}$ for $d<\beta<\alpha d$ and will play an important role in the small-balls scaling. 
Recall the following properties on ${\cal M}_{\alpha,\beta}$ from Propositions~2.2 and~2.3 from~\cite{BD2009}: 
%cf. preuve de la Prop. 2.3
\begin{prop}
\label{prop:Mcal}
\begin{enumerate}[(i)]
\item The set ${\cal M}_{\alpha,\beta}$ is a linear subspace of $\mathcal{Z}(\rd)$ and, for all $\mu\in {\cal M}_{\alpha,\beta}$,
\begin{equation*}
\label{maj:chap3}
 \int_{\real^d\times\real_+} \vert\mu\big(B(x,r)\big)\vert^\alpha r^{-\beta-1}\ dxdr <+\infty.
\end{equation*} 
\item If $d<\beta<\alpha d$, then $L^1(\real^d)\cap L^\alpha(\real^d)\subset {\cal M}_{\alpha,\beta}$ and for all $\mu\in L^1(\real^d)\cap L^\alpha (\real^d) $:
\begin{equation*}
\label{prop:L1L2}
\int_{\real^d}\vert\mu\big(B(x,r)\big)\vert^\alpha\ dx\leq C_\mu\big( r^{d}\wedge r^{\alpha  d}\big).
\end{equation*}
\end{enumerate}
\end{prop}
In the sequel, we investigate the behaviour of \eqref{expr:fluct} when $\rho \longrightarrow 0$ in various situations.

%%%%%%%%%%%%%%%%%%%%%%%%%%%%%%%%%%%%%%%%%%%%%%%%%%%%%%%%%%%%%%%%%%%%%%%%%%%%%%%%%%%%%%%%%%%%%%%%%%%%%
%%%%%%%%%%%%%%%%%%%%%%%%%%%%%%%%%%%%%%%%%%%%%%%%%%%%%%%%%%%%%%%%%%%%%%%%%%%%%%%%%%%%%%%%%%%%%%%%%%%%%%%%%

\section{Asymptotic results}\label{sec:asympt}

In order to investigate the limit of \eqref{laplace-shot-noise} when $\rho \to 0$, we consider two main scaling scenarios: the so called local and global scenarios.\\\\
\noindent In the local scenario, the scaling is properly balanced by adjusting the model parameter in the local structure of the model, i.e. in each cluster $\lambda(\rho)\rightarrow +\infty$. In this context, the key quantity driving the different regimes appears to be the mean number $\lambda(\rho)\rho^\beta$ of large balls in each cluster (see \eqref{nb-grosse-boule}).
\\\\
In this local scaling, we scale both the radii ($\rho\rightarrow 0$) and the number of balls in each cluster (i.e. $\lambda(\rho)\rightarrow +\infty$), but not scale the mean number of clusters (i.e. $\kappa(\rho)=1$) (see Section \ref{sous-sec-local}). Heuristically, each cluster $X_c$ is a Poissonian ball model as in \cite{KLNS2007} and the asymptotics are obtained as in Theorem 2 in \cite{KLNS2007} in each such model. The whole limit of the Cox model is then a mixture of the limit random fields obtained in each Poissonian cluster $X_c$ hence the randomized Poissonian limit obtained (see Theorems \ref{th:cas1-infini}, \ref{th:cox-sc1-cas2} and \ref{th:cox-sc1-cas3}).\\\\
Next we consider a global scenario where the scaling is now balanced by adjusting the model parameter of the global cluster structure of the Cox model $\kappa (\rho)\rightarrow +\infty$. In this context, the key quantity driving the different regimes is the global mean number of large balls $\kappa(\rho)\lambda(\rho)\rho^\beta$ (see \eqref{nb-grosse-boule}). In this case, we will assume that $\kappa(\rho)\rightarrow +\infty$ but we will not suppose necessarily that $\lambda(\rho)\rightarrow +\infty$. According to the behaviour of $\kappa(\rho)\lambda(\rho)\rho^\beta$, three different regimes appear, and the cluster structure is erased at the limit $\rho\rightarrow 0$ (see Section \ref{sous-sec-global}).\\\\
A natural first step for these different scenarios is to specify the behaviour of the inner integral in \eqref{laplace-shot-noise}:
\begin{equation}\label{int:base}
\int_{\rdr+} \psi_G\left(\frac{\theta \mu(B(x,r))}{n(\rho)}\right)k_\rho(x,y) dxf\left(\frac{r}{\rho}\right) \frac{dr}{\rho}.
\end{equation}
We study \eqref{int:base} in the three different normalization settings considered: 
\begin{equation*}
(1) : n(\rho)\to +\infty,~~~ (2) : n(\rho)=1, ~~~(3) : n(\rho)\to 0.
\end{equation*}
Note that in the third point of the following Proposition \ref{prop:comportement-general-cox}, we consider smooth measures $\mu\in~L^1(\rd)\cap~L^\alpha(\rd)$. Heuristically, in the third point, large balls will disappear at the limit and it will just remain small balls, that can be seen as points on the space. To identify the limit, we need more precision on the measure, and so we consider measures with intensity in $L^1(\rd)\cap~L^\alpha(\rd)$. 
In the sequel, we shall use the technical result from \cite{BD2009} (Lemma 3.1 therein):

\begin{lemme}\label{lem:carac-stable}
Suppose $X$ is in the domain of attraction of an $\alpha$-stable law $S_\alpha(\sigma,b,0)$ for some $\alpha>1$. Then
\[\ee\left[\psi(\theta X) \right]\underset{0}{\sim}-\sigma^\alpha \vert\theta\vert^\alpha \bigg(1-i\epsilon(\theta)\tan(\pi \alpha/2)b\bigg),\]
where $\epsilon(\theta)=1$ if $\theta>0$, $\epsilon(\theta)=-1$ if $\theta<0$ and $\epsilon(0)=0$.\\
Furthermore, there is some $K>0$ such that for any $\theta\in\real$,
\[\big|\ee\left[\psi(\theta X) \right]\big| \leq K\vert\theta\vert^{\alpha}. \]
\end{lemme}
With our previous notations, we have $\psi_G(\theta)=\ee\left[\psi(\theta X) \right]$. The behaviour of \eqref{int:base} is given in the following proposition.

\begin{prop}\label{prop:comportement-general-cox}
Assume \eqref{k-bornee}, \eqref{hyp:f} and \eqref{hyp-krho} hold true.
\begin{enumerate}
\item If $n(\rho)\to +\infty$, then for $\mu\in {\cal M}_{\alpha,\beta}$, we have:
\begin{align*}
\int_{\rdr+} \psi_G\left(\frac{\theta \mu(B(x,r))}{n(\rho)}\right)k_\rho(x,y) dxf\left(\frac{r}{\rho}\right) \frac{dr}{\rho}\\
&\hspace*{-8cm}\underset{\rho\to 0}{\sim}-\sigma^\alpha\vert\theta\vert^\alpha\frac{\lambda(\rho)\rho^\beta}{n(\rho)^\alpha}\int_{\rdr+} \vert\mu(B(x,r))\vert^\alpha \bigg(1-i\epsilon(\theta \mu(B(x,r))\tan(\pi \alpha/2)b \bigg) \\
&\hspace*{1.6cm} \times k(x,y) dx C_\beta r^{-\beta -1} dr.
\end{align*}
\item If $n(\rho)=1$, then for $\mu\in {\cal M}_{\alpha,\beta}$, we have:
\begin{align*}
&\int_{\rdr+} \psi_G\bigg(\theta \mu(B(x,r))\bigg)k_\rho(x,y) dxf\left(\frac{r}{\rho}\right) \frac{dr}{\rho}\\
&\hspace*{3cm}\underset{\rho\to 0}{\sim}\lambda(\rho)\rho^\beta \int_{\rdr+} \psi_G\bigg(\theta \mu(B(x,r))\bigg)k(x,y) dx C_\beta r^{-\beta -1} dr.
\end{align*}
\item If $n(\rho)\to 0$ and $\rho^d=o(n(\rho))$, then for $\mu\in L^1(\rd)\cap L^\alpha(\rd)$, we have:
\begin{align*}
&\int_{\rdr+} \psi_G\left(\frac{\theta \mu(B(x,r))}{n(\rho)}\right)k_\rho(x,y) dxf\left(\frac{r}{\rho}\right) \frac{dr}{\rho}\\
&\hspace*{3cm}\underset{\rho\to 0}{\sim}\frac{\lambda(\rho)\rho^\beta}{n(\rho)^\gamma}\int_{\rdr+} \psi_G\left(\theta \varphi(x)r^d\right)k(x,y) dx v_d^{\gamma}C_\beta r^{-\beta -1} dr.
\end{align*}
where $\mu(dx)=\varphi(x)dx$ and $\gamma=\beta/d \in ]1,\alpha[$.
\end{enumerate}
\end{prop}
\Proof~ These results are proved with the dominated convergence theorem. In the sequel, we consider a fixed $y$ in $\rd$.\\
$1.$ Let $n(\rho)\to +\infty$ and $\mu\in {\cal M}_{\alpha,\beta}$. To prove the first equivalence, we show that:
\begin{align}\label{ref1}
&\underset{\rho\to 0}{\lim}~\frac{n(\rho)^\alpha}{\lambda(\rho)\rho^\beta}\int_{\rdr+} \psi_G\left(\frac{\theta \mu(B(x,r))}{n(\rho)}\right)k_\rho(x,y) dxf\left(\frac{r}{\rho}\right) \frac{dr}{\rho}\nonumber\\
&\hspace*{-0cm}=-\sigma^\alpha\vert\theta\vert^\alpha\int_{\rdr+} \left|\mu(B(x,r))\right|^\alpha \bigg(1-i\epsilon(\theta \mu(B(x,r))\tan(\pi \alpha/2)b \bigg) k(x,y) dx C_\beta r^{-\beta -1} dr.
\end{align}
\begin{enumerate}[(i)]
\item \textit{Convergence}: Thanks to Lemma \ref{lem:carac-stable}, we have immediately together with \eqref{hyp:f} and \eqref{hyp-krho}:
\begin{align*}
&\underset{\rho\to 0}{\lim}~\frac{n(\rho)^\alpha}{\lambda(\rho)\rho^\beta}\psi_G\left(\frac{\theta \mu(B(x,r))}{n(\rho)}\right)k_\rho(x,y) f\left(\frac{r}{\rho}\right) \frac{1}{\rho}\\
&\hspace*{1cm}=-\sigma^\alpha\vert\theta\vert^\alpha\vert\mu(B(x,r))\vert^\alpha \bigg(1-i\epsilon(\theta \mu(B(x,r))\tan(\pi \alpha/2)b \bigg)  k(x,y)C_\beta r^{-\beta -1}.
\end{align*}
\item \textit{Domination}: Using Lemma \ref{lem:carac-stable} and \eqref{hyp:f}, \eqref{hyp-krho} we have:
\begin{align}\label{refA}
&\left|\frac{n(\rho)^\alpha}{\lambda(\rho)\rho^\beta}\psi_G\left(\frac{\theta \mu(B(x,r))}{n(\rho)}\right)k_\rho(x,y) f\left(\frac{r}{\rho}\right) \frac{1}{\rho}\right|\nonumber\\
&\hspace*{7cm}\leq K \vert\theta\vert^\alpha\vert\mu(B(x,r))\vert^\alpha C_1 k(x,y)C_0 r^{-\beta -1}\nonumber\\
& \hspace*{7cm}\leq K \vert\theta\vert^\alpha\vert\mu(B(x,r))\vert^\alpha  C_1 \|k\|_{\infty}C_0 r^{-\beta -1}
\end{align}
which is independent of the parameter $\rho$ and integrable on $\rdr+$ when $\mu\in {\cal M}_{\alpha,\beta}$ thanks to Proposition \ref{prop:Mcal}. 
\end{enumerate}
The dominated convergence theorem then applies and \eqref{ref1} is obtained.\\
\noindent$2.$ Let $n(\rho)=1$ and $\mu\in {\cal M}_{\alpha,\beta}$. To prove the second equivalence, we show that:
\begin{align}\label{ref2}
&\underset{\rho\to 0}{\lim}~\frac{1}{\lambda(\rho)\rho^\beta}\int_{\rdr+} \psi_G\bigg(\theta \mu(B(x,r))\bigg)k_\rho(x,y) dxf\left(\frac{r}{\rho}\right) \frac{dr}{\rho}\nonumber\\
&\hspace*{5.5cm}=\int_{\rdr+} \psi_G\bigg(\theta \mu(B(x,r))\bigg)k(x,y) dx C_\beta r^{-\beta -1} dr.
\end{align}
\begin{enumerate}[(i)]
\item \textit{Convergence}: With \eqref{hyp:f} and \eqref{hyp-krho}, we have:
\begin{align*}
\hspace*{-0.5cm}\underset{\rho\to 0}{\lim}~\frac{1}{\lambda(\rho)\rho^\beta}\psi_G\bigg(\theta \mu(B(x,r))\bigg)k_\rho(x,y) f\left(\frac{r}{\rho}\right) \frac{1}{\rho}=\psi_G\bigg(\theta \mu(B(x,r))\bigg)k(x,y) C_\beta r^{-\beta -1}.
\end{align*}
\item \textit{Domination}: Using Lemma \ref{lem:carac-stable} and \eqref{hyp:f}, \eqref{hyp-krho}, we have:
\begin{align}\label{refB}
\hspace*{-0.5cm}\left|\frac{1}{\lambda(\rho)\rho^\beta}\psi_G\left(\theta \mu(B(x,r))\right)k_\rho(x,y) f\left(\frac{r}{\rho}\right) \frac{1}{\rho}\right|&\leq K \vert\theta\vert^\alpha\vert\mu(B(x,r))\vert^\alpha  C_1 k(x,y)C_0 r^{-\beta -1}\nonumber\\
&\leq  K \vert\theta\vert^\alpha\vert\mu(B(x,r))\vert^\alpha C_1 \|k\|_{\infty}C_0 r^{-\beta -1}
\end{align}
which is independent of the parameter $\rho$ and integrable on $\rdr+$ when $\mu\in {\cal M}_{\alpha,\beta}$ thanks to Proposition \ref{prop:Mcal}.
\end{enumerate}
 The dominated convergence theorem applies again and we obtain \eqref{ref2}.\\
\noindent $3.$ Let $n(\rho)\to 0$ and $\rho^d=o(n(\rho))$. Let $\mu\in L^1(\rd)\cap L^\alpha(\rd)$ with $\mu(dx)=\varphi(x)dx$. To prove the last equivalence, we show that:
\begin{align}\label{ref3}
&\underset{\rho\to 0}{\lim}~\frac{n(\rho)^\gamma}{\lambda(\rho)\rho^\beta}\int_{\rdr+} \psi_G\left(\frac{\theta \mu(B(x,r))}{n(\rho)}\right)k_\rho(x,y) dxf\left(\frac{r}{\rho}\right) \frac{dr}{\rho}\nonumber\\
&\hspace*{6.5cm}=\int_{\rdr+} \psi_G\left(\theta \varphi(x)r^d\right)k(x,y) dx v_d^{\gamma}C_\beta r^{-\beta -1} dr.
\end{align}
First, the change of variable $r=n(\rho)^{1/d}s$ gives:
\begin{align*}
&\frac{n(\rho)^\gamma}{\lambda(\rho)\rho^\beta}\int_{\rdr+} \psi_G\left(\frac{\theta \mu(B(x,r))}{n(\rho)}\right)k_\rho(x,y) dxf\left(\frac{r}{\rho}\right) \frac{dr}{\rho}\\
&\hspace*{1cm}=\frac{n(\rho)^\gamma}{\lambda(\rho)\rho^\beta}\int_{\rdr+} \psi_G\left(\frac{\theta \mu(B(x,n(\rho)^{1/d}s))}{n(\rho)}\right)k_\rho(x,y) dxf\left(\frac{n(\rho)^{1/d}s}{\rho}\right) \frac{n(\rho)^{1/d}}{\rho}ds,
\end{align*}
and we study the limit of this latter expression like in 1 and 2:
\begin{enumerate}[(i)]
\item \textit{Convergence}:  First, since $\varphi\in L^1(\rd)$, Lemma 4 in \cite{KLNS2007} applies and for all $y\in\rd$ we have:
\begin{align*}
&\underset{\rho\rightarrow 0}{\lim}~\frac{n(\rho)^\gamma}{\lambda(\rho)\rho^\beta} \psi_G\left(\frac{\theta \mu(B(x,n(\rho)^{1/d}s))}{n(\rho)}\right)k_\rho(x,y) f\left(\frac{n(\rho)^{1/d}s}{\rho}\right) \frac{n(\rho)^{1/d}}{\rho}\\
&\hspace*{8cm}=\psi_G\left(\theta \varphi(x)v_d s^d\right)k(x,y)C_\beta s^{-\beta -1}.
\end{align*}
where $v_d$ is the Lebesgue measure of the unit ball of $\rd$. 
\item \textit{Domination}: Like in \cite{KLNS2007}, set $\displaystyle\varphi_*(x)=\underset{v>0}{\sup}~\frac{\mu(B(x,v))}{v^d}$. \\
Thanks to Lemma \ref{lem:carac-stable}, we have $\left|\psi_G(u)\right|\leq K\vert u\vert^{\alpha}$ and combining with \eqref{ineg-psig}, we have
\begin{equation}\label{eq:ineg psig 1}
\left|\psi_G(u)\right|\leq M \left(\vert u\vert \wedge\vert u\vert^{\alpha}\right)
\end{equation}
for some constant $M>0$. \\
Let $\varepsilon >0$ such that $1< \gamma-\varepsilon<\gamma+\varepsilon<\alpha$. Thanks to \eqref{eq:ineg psig 1}, we have
\begin{equation}\label{eq:ineg psig 2}
\left|\psi_G(u)\right|\leq M \left(\vert u\vert^{\gamma-\varepsilon} \wedge\vert u\vert^{{\gamma+\varepsilon} }\right).
\end{equation}
Since $\varphi\in L^1(\rd)\cap L^\alpha(\rd)$, $\varphi \in L^q(\rd)$ for all $1\leq q\leq \alpha$ so that Lemma 4 in \cite{KLNS2007} ensures $\displaystyle\varphi_*\in L^{q}(\rd)$ for all $1<q<\alpha$ which entails $\displaystyle\varphi_*\in L^{\gamma-\varepsilon}(\rd)\cap L^{\gamma+\varepsilon}(\rd)$.\\
Thanks to \eqref{k-bornee}, \eqref{hyp:f}, \eqref{hyp-krho} and \eqref{eq:ineg psig 2}, for $\rho>0$ and $\theta \in \real$ we have the following domination:
\begin{align}\label{eq:controle-petite-boule}
&\hspace{-1cm}\left|\frac{n(\rho)^\gamma}{\lambda(\rho)\rho^\beta}\psi_G\left(\frac{\theta \mu(B(x,n(\rho)^{1/d}s))}{n(\rho)}\right)k_\rho(x,y)f\left(\frac{n(\rho)^{1/d}s}{\rho}\right) \frac{n(\rho)^{1/d}}{\rho}\right|\nonumber \\
&\hspace{-1cm}\leq M\frac{n(\rho)^\gamma}{\lambda(\rho)\rho^\beta}\left\lbrace\left|\frac{\theta \mu(B(x,n(\rho)^{1/d}s))}{n(\rho)}\right|^{\gamma-\varepsilon}\wedge \left|\frac{\theta \mu(B(x,n(\rho)^{1/d}s))}{n(\rho)}\right|^{\gamma+\varepsilon}\right\rbrace \nonumber\\
&\hspace{8.3cm}\times C_1 \lambda(\rho) \|k\|_{\infty}C_0 \frac{\rho^{\beta}}{n(\rho)^{\gamma}} s^{-\beta -1} \nonumber \\
&\hspace{-1cm}\leq M C_1 C_0 \|k\|_{\infty}\left|\theta  \varphi_*(x) s^d\right|^{\gamma-\varepsilon}\wedge \left|\theta  \varphi_*(x) s^d\right|^{\gamma+\varepsilon}  s^{-\beta -1} \nonumber\\
&\hspace{-1cm}\leq M \left(\vert\theta\vert ^{\gamma-\varepsilon}+\vert\theta\vert ^{\gamma+\varepsilon}\right)C_1 C_0 \|k\|_{\infty} \left(\vert\varphi_*(x)\vert^{\gamma-\varepsilon}+ \vert\varphi_*(x)\vert^{\gamma+\varepsilon}\right)\left(s^{-\varepsilon d -1}\wedge s^{\varepsilon d -1}  \right)
\end{align}
\hspace{-1.1cm} which is independent of the parameter $\rho >0$ and integrable on $\rdr+$. 
\end{enumerate}
Finally, the dominated convergence theorem applies one more time and we obtain:
\begin{align*}
&\underset{\rho\to 0}{\lim}~\frac{n(\rho)^\gamma}{\lambda(\rho)\rho^\beta}\int_{\rdr+} \psi_G\left(\frac{\theta \mu(B(x,r))}{n(\rho)}\right)k_\rho(x,y) dxf\left(\frac{r}{\rho}\right) \frac{dr}{\rho}\\
&\hspace*{6.5cm}=\int_{\rdr+} \psi_G\left(\theta \varphi(x)v_d s^d\right)k(x,y)dxC_\beta s^{-\beta -1}ds
\end{align*}
which gives \eqref{ref3}, up to the change of variable $r=v_d^{1/d}s$.
\CQFD
\\\\
In the sequel, we establish convergences of finite-dimensional distributions by proving one-dimensional convergences of the distributions and using the Cram\'er-Wold  device combined with the linear structure of ${\cal M}_{\alpha,\beta}$.
%%%%%%%%%%%%%%%%%%%%%%%%%%%%%%%%%%%%%%%%%%%%%%%%%%%%%%%%%%%%%%%%%%%%%%%%%%%%%%%%%%%%%%%%%%%%%%%%%%%%%

\subsection{Local scaling}\label{sous-sec-local}

In this section, we investigate the so called local scenario where the radii together with the mean number of large balls in each cluster is rescaled. However, here, the mean number of clusters remains constant, i.e. $\kappa(\rho)=1$, and in this case the intensity of the shot-noise Cox process specifies as
\begin{equation*}
 Z_\rho(x)=\sum_{y\in\Phi}k_\rho(x,y),
\end{equation*} 
where $\Phi$ is a Poisson point process on $\rd$ with intensity measure $dy$.
We assume that $\underset{\rho\to 0}{\lim}~\lambda (\rho)=+\infty$ i.e. while we zoom-out ($\rho\rightarrow 0$), the clusters are bigger and bigger ($\lambda (\rho)\rightarrow +\infty$) but the mean number of clusters does not change.

\begin{theo}\label{th:cas1-infini}
Assume \eqref{k-bornee}, \eqref{eq:k-intensity}, \eqref{hyp:f} and \eqref{hyp-krho} hold true. Suppose $\lambda(\rho)\rho^\beta \underset{\rho \rightarrow 0}{\longrightarrow}+\infty$ and set $n(\rho)=\left(\lambda(\rho)\rho^\beta\right)^{1/\alpha}$. Then the following limit holds when $\rho\to 0$:
\begin{equation}\label{th:sc1-LB}
\frac{M_\rho(\mu)-\ee\left[M_\rho(\mu)\left.\right|\Lambda_\rho\right]}{n(\rho)}\underset{fdd}{\overset{{\cal M}_{\alpha,\beta}}{\Longrightarrow}} \int_{\rdr+}\mu(B(x,r))M_\alpha(dx,dr)
\end{equation}
where conditionally to $L_\alpha$, $M_\alpha$ is a $\alpha$-stable random measure with control measure $L_\alpha(dx,dr)=\sigma^\alpha Z(x)dx C_\beta r^{-\beta -1}dr$, where $Z$ is given in \eqref{def:Z}-\eqref{eq:k-intensity}, and constant skewness function $b$.
\end{theo}
\Proof~ We apply the dominated convergence theorem to take the limit when $\rho\to 0$ in \eqref{laplace-shot-noise}. From (1) in Proposition \ref{prop:comportement-general-cox}, we have:
\begin{align*}
\underset{\rho\rightarrow 0}{\lim}~~1-&\exp\left(\int_{\rdr+} \psi_G\left(\frac{\theta \mu(B(x,r))}{n(\rho)}\right)k_\rho(x,y) dxf\left(\frac{r}{\rho}\right) \frac{dr}{\rho}\right) \nonumber\\
&\hspace{-1.5cm}=1-\exp\bigg(-\sigma^\alpha\vert\theta\vert^\alpha\int_{\rdr+} \vert\mu(B(x,r))\vert^\alpha \big(1-i\epsilon(\theta \mu(B(x,r))\tan(\pi \alpha/2)b \big)\\
&\hspace*{9cm}\times k(x,y) dx C_\beta r^{-\beta -1} dr\bigg).
\end{align*}
In order to prove the domination, let $h(u)=1-e^u$, $u\in\mathbb{C}$. By the mean value theorem, we have $\left|h(u)\right|=\left|e^u-1\right|\leq e^A \vert u\vert$ for all $\vert  u\vert\leq A$, where $A$ is a fixed positive constant. In our context we take:
$$\displaystyle u=\int_{\rdr+} \psi_G\left(\frac{\theta \mu(B(x,r))}{n(\rho)}\right)k_\rho(x,y) dxf\left(\frac{r}{\rho}\right) \frac{dr}{\rho}$$
and
$$\displaystyle A:=K\|k\|_{\infty}C_0C_1 \vert\theta\vert^\alpha  \int_{\rdr+}\vert\mu(B(x,r))\vert^\alpha r^{-\beta -1}dxdr<+\infty.$$
Using the same domination as \eqref{refA} in the proof of (1) in Proposition \ref{prop:comportement-general-cox}, we have indeed $\vert u\vert \leq A$ and consequently:
\begin{align*}
&\left|1-\exp\left(\int_{\rdr+} \psi_G\left(\frac{\theta \mu(B(x,r))}{n(\rho)}\right)k_\rho(x,y) dxf\left(\frac{r}{\rho}\right) \frac{dr}{\rho}\right)\right| \nonumber\\
&\hspace{4cm}\leq e^{A}\left|\int_{\rdr+} \psi_G\left(\frac{\theta \mu(B(x,r))}{n(\rho)}\right)k_\rho(x,y) dxf\left(\frac{r}{\rho}\right) \frac{dr}{\rho}\right| \nonumber\\
&\hspace{4cm}\leq e^{A}K C_0 C_1\vert\theta\vert ^\alpha\int_{\rdr+} \vert\mu(B(x,r))\vert^\alpha k(x,y) r^{-\beta -1}dxdr 
\end{align*}
again with Lemma \ref{lem:carac-stable}. Since the bound is independent of the parameter $\rho$ and is integrable on $\rd$ with respect to the Lebesgue measure thanks to Fubini theorem, condition \eqref{eq:k-intensity} and Proposition \ref{prop:Mcal}, the dominated convergence theorem applies.\\\\
As a consequence, under the condition of Theorem \ref{th:cas1-infini} the limit in \eqref{laplace-shot-noise} writes:
\begin{align*}
&\underset{\rho\rightarrow 0}{\lim}\ee\left[\exp\left(-\theta\left(\frac{M_\rho(\mu)-\ee\left[M_\rho(\mu)\left.\right|\Lambda_\rho\right]}{n(\rho)}  \right)\right) \right]\nonumber \\
&\hspace{2.5cm}=\exp\bigg(-\int_{\rd} \bigg( 1-\exp\bigg(-\sigma^\alpha\vert\theta\vert^\alpha\int_{\rdr+} \vert\mu(B(x,r))\vert^\alpha \\
&\hspace*{4.5cm}\times(1-i\epsilon(\theta \mu(B(x,r))\tan(\pi \alpha/2)b ) k(x,y) dx C_\beta r^{-\beta -1} dr\bigg) \bigg)dy \bigg).
\end{align*}
We identify the obtained limit as the characteristic function of $$\displaystyle\int_{\rdr+}\mu(B(x,r))M_\alpha(dx,dr)$$ where $M_\alpha$ is as described in Theorem \ref{th:cas1-infini} (we refer \cite{ST1994} for basics on stable measure and stable random variable). Indeed, for $\theta \in\real$ and $\mu\in {\cal M}_{\alpha,\beta}$, we have:
\begin{align*}
&\ee\left[\exp\left(i\theta \int_{\rdr+}\mu(B(x,r))M_\alpha(dx,dr) \right) \right]\\
&\hspace*{0cm}=\ee\left[ \ee\left[\left.\exp\left(i\theta \int_{\rdr+}\mu(B(x,r))M_\alpha(dx,dr) \right) \right| L_\alpha\right] \right]\\
&\hspace*{0cm}=\ee\bigg[\exp\bigg(  - \sigma^\alpha\vert\theta\vert^\alpha\int_{\rdr+} \vert\mu(B(x,r))\vert^\alpha \left(1-i\epsilon(\theta \mu(B(x,r))\tan(\pi \alpha/2)b \right)\\
&\hspace*{10.5cm}\times Z(x) dx C_\beta r^{-\beta -1} dr \bigg) \bigg]\\
&\hspace*{0cm}=\ee\bigg[\exp\bigg(\int_{\rd} \bigg(  \int_{\rdr+}  - \sigma^\alpha\vert\theta\vert^\alpha    \vert\mu(B(x,r))\vert^\alpha \left(1-i\epsilon(\theta \mu(B(x,r))\tan(\pi \alpha/2)b \right) \\
&\hspace*{9cm}\times k(x,y) dx C_\beta r^{-\beta -1} dr \bigg)\Phi(dy) \bigg) \bigg]\\
&\hspace*{0cm}=\exp\bigg(-\int_{\rd} \bigg( 1-\exp\bigg(-\sigma^\alpha\vert\theta\vert^\alpha\int_{\rdr+} \vert\mu(B(x,r))\vert^\alpha \\
&\hspace*{4.5cm}\times(1-i\epsilon(\theta \mu(B(x,r))\tan(\pi \alpha/2)b ) k(x,y) dx C_\beta r^{-\beta -1} dr\bigg) \bigg)dy \bigg)
\end{align*}
by Proposition \ref{prop:laplace-poisson} since $\Phi$ is a Poisson point process on $\rd$ with intensity measure $dy$.\\
Finally, we have:
\begin{align*}
&\underset{\rho\rightarrow 0}{\lim}\ee\left[\exp\left(i\theta\left(\frac{M_\rho(\mu)-\ee\left[M_\rho(\mu)\left.\right|\Lambda_\rho\right]}{n(\rho)}  \right)\right) \right]\\
&\hspace*{6.5cm}=\ee\left[\exp\left(i\theta \int_{\rdr+}\mu(B(x,r))M_\alpha(dx,dr) \right) \right]
\end{align*}
which proves \eqref{th:sc1-LB}.
\CQFD

\begin{theo}\label{th:cox-sc1-cas2}
Assume \eqref{k-bornee}, \eqref{eq:k-intensity}, \eqref{hyp:f} and \eqref{hyp-krho} hold true.\\ Suppose $\lambda(\rho)\rho^\beta \underset{\rho \rightarrow 0}{\longrightarrow}a \in ]0,+\infty[$ and set $n(\rho)=1$. Then the following limit holds when $\rho\to 0$:
\begin{equation*}
M_\rho(\mu)-\ee\left[M_\rho(\mu)\left.\right|\Lambda_\rho\right]\underset{fdd}{\overset{{\cal M}_{\alpha,\beta}}{\Longrightarrow}}N_\rho(\mu)-\ee\left[N_\rho(\mu)\left.\right|\Lambda'_\rho\right]
\end{equation*}
where $\displaystyle N_\rho(\mu)=\int_{\rrr}m\mu(B(x,r))C'(dx,dr,dm)$, $C'$ is a Cox process on $\rrr$ directed by $\Lambda'(dx,dr,dm)=Z(x)dxaC_\beta r^{-\beta-1}drG(dm)$, and $Z$ is given in \eqref{def:Z}-\eqref{eq:k-intensity}.
\end{theo}
\Proof~ Using (2) in Proposition \ref{prop:comportement-general-cox}, for all $y\in\rd$, we have:
\begin{align*}
\underset{\rho\rightarrow 0}{\lim}~~1-&\exp\left(\int_{\rdr+} \psi_G\bigg(\theta \mu(B(x,r))\bigg)k_\rho(x,y) dxf\left(\frac{r}{\rho}\right) \frac{dr}{\rho}\right) \nonumber\\
&\hspace{3cm}=1-\exp\left(\int_{\rdr+} \psi_G\bigg(\theta \mu(B(x,r))\bigg)k(x,y) dxaC_\beta r^{-\beta -1}dr\right). 
\end{align*}
Consider like in the proof of Theorem \ref{th:cas1-infini} $h(u)=1-e^u$, $u\in\mathbb{C}$. Let $\rho_1>0$ such that $\lambda(\rho)\rho^\beta \leq2a$ for all $0<\rho<\rho_1$. In our context, take:
$$ u=\int_{\rdr+} \psi_G\bigg(\theta \mu(B(x,r))\bigg)k_\rho(x,y) dxf\left(\frac{r}{\rho}\right) \frac{dr}{\rho}$$
and set
$$\displaystyle A:=a K \|k\|_{\infty}C_0 C_1 \vert\theta\vert^\alpha  \int_{\rdr+}\vert\mu(B(x,r))\vert^\alpha r^{-\beta -1}dxdr<+\infty.$$
Using the very same domination \eqref{refB} as in the proof of (2) in Proposition \ref{prop:comportement-general-cox}, we have indeed $\vert u\vert\leq A$ when $0<\rho<\rho_1$. Thereby, for all $0<\rho<\rho_1$:
\begin{align*}
&\left|1-\exp\left(\int_{\rdr+} \psi_G\left(\theta \mu(B(x,r))\right)k_\rho(x,y) dxf\left(\frac{r}{\rho}\right) \frac{dr}{\rho}\right)\right| \nonumber\\
&\hspace{4cm}\leq e^A\left|\int_{\rdr+} \psi_G\left(\theta \mu(B(x,r))\right)k_\rho(x,y) dxf\left(\frac{r}{\rho}\right) \frac{dr}{\rho}\right| \nonumber\\
&\hspace{4cm}\leq e^{A}a K C_0 C_1 \vert\theta\vert^\alpha\int_{\rdr+} \vert\mu(B(x,r))\vert^\alpha k(x,y) r^{-\beta -1}dxdr .
\end{align*}
Using again Lemma \ref{lem:carac-stable}, since the bound is independent of the parameter $\rho\in ]0,\rho_1[$ and is integrable on $\rd$ with respect to the Lebesgue measure thanks to Fubini, \eqref{eq:k-intensity} and Proposition~\ref{prop:Mcal}, the dominated convergence theorem applies and the limit in \eqref{laplace-shot-noise} writes
\begin{align*}
\underset{\rho\rightarrow 0}{\lim}&\ee\bigg[\exp\bigg(i\theta\left(M_\rho(\mu)-\ee\left[M_\rho(\mu)\left.\right|\Lambda_\rho\right] \right)\bigg) \bigg]\nonumber \\
&\hspace{0cm}=\exp\left(-\int_{\rd} \left(1-\exp\left(\int_{\rdr+} \psi_G\left(\theta \mu(B(x,r))\right)k(x,y) dxaC_\beta r^{-\beta -1}dr\right) \right)dy \right)
\end{align*}
which is the characteristic function of $\displaystyle N_\rho(\mu)-\ee\left[N_\rho(\mu)\left.\right|\Lambda'_\rho\right]$ (see Proposition \ref{prop:laplace-cox-centree}). This proves Theorem \ref{th:cox-sc1-cas2}. \CQFD

\begin{theo}\label{th:cox-sc1-cas3}
Assume \eqref{k-bornee}, \eqref{eq:k-intensity}, \eqref{hyp:f} and \eqref{hyp-krho} hold true.\\ Suppose $\lambda(\rho)\rho^\beta \underset{\rho \rightarrow 0}{\longrightarrow}0$ and set $n(\rho)=\left(\lambda(\rho)\rho^\beta\right)^{1/\gamma}$, where $\gamma=\beta /d \in ]1,\alpha[$. Then the following limit holds when $\rho\to 0$:
\begin{equation*}
\frac{M_\rho(\mu)-\ee\left[M_\rho(\mu)\left.\right|\Lambda_\rho\right]}{n(\rho)}\underset{fdd}{\overset{L^1(\rd)\cap L^\alpha(\rd)}{\Longrightarrow}}  \int_{\rd}\varphi (x)M_\gamma (dx),
\end{equation*}
for $\mu(dx)=\varphi(x)dx$, where conditionally to $S$, $M_\gamma$ is a $\gamma$-stable measure with control measure $S(dx)=\sigma_\gamma Z(x)dx$ for
$$\displaystyle\sigma_\gamma=\frac{C_\beta v_d^\gamma}{d}\bigg(\int_0^{+\infty}\frac{1-\cos(r)}{r^{1+\gamma}}dr\bigg)\bigg(\int_\real \vert m\vert^\gamma G(dm)\bigg)$$ 
with constant skewness function equals to
$$b_\gamma=-\frac{\int_\real \epsilon (m)\vert m\vert^\gamma G(dm)}{\int_\real \vert m\vert^\gamma G(dm)}$$
and $Z$ given in \eqref{def:Z}-\eqref{eq:k-intensity}.
\end{theo}
\Proof~ We apply the dominated convergence theorem. Using (3) in Proposition \ref{prop:comportement-general-cox}, for all $y\in\rd$, we have:
\begin{align*}
\underset{\rho\rightarrow 0}{\lim}~~1-&\exp\left(\int_{\rdr+} \psi_G\left(\frac{\theta \mu(B(x,r))}{n(\rho)}\right)k_\rho(x,y) dxf\left(\frac{r}{\rho}\right) \frac{dr}{\rho}\right) \nonumber\\
&\hspace{3cm}=1-\exp\left(\int_{\rdr+} \psi_G\left(\theta \varphi(x)r^d\right)k(x,y) dx v_d^{\gamma}C_\beta r^{-\beta -1} dr\right). 
\end{align*}
In order to derive the domination, let $\varepsilon >0$ be such that $1< \gamma-\varepsilon<\gamma+\varepsilon<\alpha$. Consider again $h(u)=1-e^u$, $u\in\mathbb{C}$ and, in our context, take
\begin{align*}
&A=M\left(\vert\theta\vert ^{\gamma-\varepsilon}+\vert\theta\vert ^{\gamma+\varepsilon}\right)\|k\|_{\infty}C_0 C_1 \\
&\hspace*{3.5cm} \times \int_{\rdr+} \left(\varphi_*(x)^{\gamma-\varepsilon}+ \varphi_*(x)^{\gamma+\varepsilon}\right)\left(r^{-\varepsilon d -1}\wedge r^{\varepsilon d -1}\right) dxdr<+\infty,
\end{align*}
and
$$ u=\int_{\rdr+} \psi_G\left(\frac{\theta \mu(B(x,r))}{n(\rho)}\right)k_\rho(x,y) dxf\left(\frac{r}{\rho}\right) \frac{dr}{\rho}.$$
Using the very same domination \eqref{eq:controle-petite-boule} as in the proof of (3) in Proposition \ref{prop:comportement-general-cox}, we have $\vert u\vert\leq A$ and for all $\rho>0$:
\begin{align*}
&\left|1-\exp\left(\int_{\rdr+} \psi_G\left(\frac{\theta \mu(B(x,r))}{n(\rho)}\right)k_\rho(x,y) dxf\left(\frac{r}{\rho}\right) \frac{dr}{\rho}\right)\right| \\
&\hspace{0cm}\leq e^{A}\left|\int_{\rdr+}  \psi\left(\frac{\theta \mu(B(x,r))}{n(\rho)}\right)k_\rho(x,y) dxf\left(\frac{r}{\rho}\right) \frac{dr}{\rho}\right| \\
&\hspace{0cm}\leq M e^{A}\left(\vert\theta\vert ^{\gamma-\varepsilon}+\vert\theta\vert ^{\gamma+\varepsilon}\right)C_0 C_1\\
&\hspace*{4cm} \times\int_{\rdr+} \left(\varphi_*(x)^{\gamma-\varepsilon}+ \varphi_*(x)^{\gamma+\varepsilon}\right)k(x,y)\left(r^{-\varepsilon d -1}\wedge r^{\varepsilon d -1}\right) dxdr. 
\end{align*}
Like for \eqref{eq:controle-petite-boule}, the bound is finite, and independent of the parameter $\rho >0$ and integrable on $\rd$ with respect to the Lebesgue measure thanks to Fubini theorem and condition \eqref{eq:k-intensity}.
The dominated convergence theorem applies and the limit in \eqref{laplace-shot-noise} writes:
\begin{align*}
\underset{\rho\rightarrow 0}{\lim}&\ee\left[\exp\left(-\theta\left(\frac{M_\rho(\mu)-\ee\left[M_\rho(\mu)\left.\right|\Lambda_\rho\right]}{n(\rho)}  \right)\right) \right]\nonumber \\
&\hspace{0cm}=\exp\left(-\int_{\rd} \left(1-\exp\left(   \int_{\rdr+}\psi_G\left(\theta \varphi(x)r^d\right)k(x,y)dxv_d^{\gamma}C_\beta r^{-\beta -1}dr   \right) \right)dy \right).
\end{align*}
We identify the obtained limit as the characteristic function of $\displaystyle \int_{\rd}\varphi (x)M_\gamma (dx)$ where $M_\gamma$ is described in Theorem \ref{th:cox-sc1-cas3}. Indeed, for $\theta\in\real$ and $\mu\in L^1(\rd)\cap L^\alpha(\rd)$, we have:
\begin{align*}
&\ee\left[\exp\left(i\theta \int_{\rd}\varphi (x)M_\gamma (dx) \right) \right]\\
&\hspace*{1.3cm}=\ee\left[ \ee\left[\left.\exp\left(i\theta \int_{\rd}\varphi (x)M_\gamma (dx) \right) \right| S\right] \right]\\
&\hspace*{1.3cm}=\ee\left[\exp\left(\int_{\rdr+}\psi_G\left(\theta \varphi(x)r^d\right)Z(x)dxv_d^{\gamma}C_\beta r^{-\beta -1}dr  \right) \right]\\
&\hspace*{1.3cm}=\ee\left[\exp\left(\int_{\rd} \left(  \int_{\rdr+}\psi_G\left(\theta \varphi(x)r^d\right)k(x,y)dxv_d^{\gamma}C_\beta r^{-\beta -1}dr  \right)\Phi(dy) \right) \right]\\
&\hspace*{1.3cm}=\exp\left(-\int_{\rd} \left(1-\exp\left(   \int_{\rdr+}\psi_G\left(\theta \varphi(x)r^d\right)k(x,y)dxv_d^{\gamma}C_\beta r^{-\beta -1}dr   \right) \right)dy \right)
\end{align*}
since $\Phi$ is a Poisson point process on $\rd$ with intensity measure $dy$.\\
Finally, we have:
\[\underset{\rho\rightarrow 0}{\lim}\ee\left[\exp\left(i\theta\left(\frac{M_\rho(\mu)-\ee\left[M_\rho(\mu)\left.\right|\Lambda_\rho\right]}{n(\rho)}  \right)\right) \right]=\ee\left[\exp\left(i\theta \int_{\rd}\varphi (x)M_\gamma (dx) \right) \right],\]
proving Theorem \ref{th:cox-sc1-cas3}.
\CQFD

\begin{Rem}
In this section, the results obtained are somehow a randomnization of the corresponding results in the Poissonian case. Actually, the results are the same but with random intensity. This is due to the considered shot-noise model that exhibits a clusters structure, and with this scenario of scaling, we do not scale the mean number of cluster. We perform a zoom-out in each cluster, and the result is heuristically the mixture of the different limits obtained, where the location of the limits are, here, random.
\end{Rem}

\subsection{Global scaling}\label{sous-sec-global}

In this section, we perform a global scaling on the model. Heuristically, as in the Poissonian model, we focus on the mean number of large balls in the whole model. In the Poissonian or determinantal case, this key quantity is $\lambda(\rho)\rho^\beta$. In the shot-noise Cox model described below, a similar computation shows that the mean number of large balls is of order $\kappa(\rho)\lambda(\rho)\rho^\beta$ (see~\eqref{nb-grosse-boule}). The behaviour of this quantity drives the fluctuations of our Cox model. \\
In this context, we assume that $\kappa(\rho)\longrightarrow +\infty$ but do not impose that $\lambda(\rho)\longrightarrow +\infty$.\\\\
In the sequel, we use the following elementary observation: if $u_\rho \underset{\rho\rightarrow 0}{\sim}v_\rho$ and $v_\rho \underset{\rho\rightarrow 0}{\longrightarrow}0$ then:
\begin{equation}\label{eq:equiv-1-exp}
1-e^{u_\rho}\underset{\rho\rightarrow 0}{\sim}-v_\rho.
\end{equation}

\begin{theo}\label{th:cas global}
Assume \eqref{k-bornee}, \eqref{eq:k-intensity}, \eqref{hyp:f} and \eqref{hyp-krho} hold true.
\begin{enumerate}
\item Suppose $\kappa(\rho)\lambda(\rho)\rho^\beta \underset{\rho \rightarrow 0}{\longrightarrow}+\infty$ and set $n(\rho)=\left(\kappa(\rho)\lambda(\rho)\rho^\beta\right)^{1/\alpha}$. Then the following limit holds when $\rho\to 0$:
\begin{equation*}
\frac{M_\rho(\mu)-\ee\left[M_\rho(\mu)\left.\right|\Lambda_\rho\right]}{n(\rho)}\underset{fdd}{\overset{{\cal M}_{\alpha,\beta}}{\Longrightarrow}} \int_{\rdr+}\mu(B(x,r))\widetilde{M_\alpha}(dx,dr)
\end{equation*}
where $\widetilde{M_\alpha}$ is a $\alpha$-stable random measure with control measure $\sigma^\alpha dx C_\beta r^{-\beta -1}dr$ and constant skewness function $b$.

\item Suppose $\kappa(\rho)\lambda(\rho)\rho^\beta \underset{\rho \rightarrow 0}{\longrightarrow}a \in]0,+\infty[$ and set  $n(\rho)=1$. Then the following limit holds when $\rho\to 0$:
\begin{equation*}
M_\rho(\mu)-\ee\left[M_\rho(\mu)\left.\right|\Lambda_\rho\right]\underset{fdd}{\overset{{\cal M}_{\alpha,\beta}}{\Longrightarrow}} \int_{\rrr} m\mu(B(x,r))\widetilde{\Pi}(dx,dr,dm)
\end{equation*}
where $\widetilde{\Pi}$ is a centered Poisson random measure with control measure \\$aC_\beta r^{-\beta -1}dxdrG(dm)$.

\item Suppose $\kappa(\rho)\lambda(\rho)\rho^\beta \underset{\rho \rightarrow 0}{\longrightarrow}0$ and $\kappa(\rho)\lambda(\rho)\longrightarrow +\infty.$\\
Set  $n(\rho)=\left(\kappa(\rho)\lambda(\rho)\rho^\beta\right)^{\frac{1}{\gamma}}$ with $\gamma=\beta/d\in]1,\alpha[$. Then the following limit holds when $\rho\to 0$:
\begin{equation*}
\frac{M_\rho(\mu)-\ee\left[M_\rho(\mu)\left.\right|\Lambda_\rho\right]}{n(\rho)}\underset{fdd}{\overset{L^1(\rd)\cap L^\alpha(\rd)}{\Longrightarrow}} \int_{\rd}\varphi (x)M_\gamma(dx)
\end{equation*}
for $\mu(dx)=\varphi(x)dx$, where $M_\gamma$ is a $\gamma$-stable measure with control measure $\sigma_\gamma dx$ and constant unit skewness $b_\gamma$ given in Theorem \ref{th:cox-sc1-cas3}.
\end{enumerate} 
\end{theo}

\begin{Rem} Results 1. and 2. are very closed to the Poissonian case in \cite{BEK2010}. It is important to note that for the third result, given that we do not necessarily have $\lambda(\rho)\longrightarrow +\infty$, we must impose the constraint $\kappa(\rho)\lambda(\rho)\longrightarrow +\infty$ that allows us to adjust the speed of increase of the cluster number according to the speed at which the number of balls in a cluster varies. Of course, if $\lambda(\rho)\longrightarrow +\infty$, this condition is necessarily verified, but if the number of balls in a cluster remains constant or goes to 0, this condition tells us how fast the number of clusters should increase to have a non trivial limit.
\end{Rem}

\noindent\Proof~ In this proof, we skip the major part of the details. The proofs of the previous theorems contain all the elements to justify the following results. Here we just give the limit when $\rho\rightarrow 0$ in \eqref{int:base} in the different cases.\\
1. From (1) in Proposition \ref{prop:comportement-general-cox}, for $\mu\in {\cal M}_{\alpha,\beta}$, we have:
\begin{align*}
&\int_{\rdr+} \psi_G\left(\frac{\theta \mu(B(x,r))}{n(\rho)}\right)k_\rho(x,y) dxf\left(\frac{r}{\rho}\right) \frac{dr}{\rho}\\
&\hspace*{0cm}\underset{\rho\to 0}{\sim}-\sigma^\alpha\vert\theta\vert^\alpha\frac{\lambda(\rho)\rho^\beta}{n(\rho)^\alpha}\int_{\rdr+} \vert\mu(B(x,r))\vert^\alpha \bigg(1-i\epsilon(\theta \mu(B(x,r))\tan(\pi \alpha/2)b \bigg)\\
&\hspace*{10.5cm}\times k(x,y) dx C_\beta r^{-\beta -1} dr\\
&=-\sigma^\alpha\vert\theta\vert^\alpha\frac{1}{\kappa(\rho)}\int_{\rdr+} \vert\mu(B(x,r))\vert^\alpha \bigg(1-i\epsilon(\theta \mu(B(x,r))\tan(\pi \alpha/2)b \bigg)\\
&\hspace*{10.5cm}\times k(x,y) dx C_\beta r^{-\beta -1} dr.
\end{align*}
Since $\displaystyle\underset{\rho\to 0}{\lim}~\kappa(\rho)=+\infty$, \eqref{eq:equiv-1-exp} gives with \eqref{hyp:f} and \eqref{hyp-krho}:
\begin{align*}
&\underset{\rho\to 0}{\lim}~\left(1-\exp\left(\int_{\rdr+} \psi_G\left(\frac{\theta \mu(B(x,r))}{n(\rho)}\right)k_\rho(x,y) dxf\left(\frac{r}{\rho}\right) \frac{dr}{\rho}\right) \right)\kappa(\rho)\\
&\hspace*{1.7cm}=\sigma^\alpha\vert\theta\vert^\alpha\int_{\rdr+} \vert\mu(B(x,r))\vert^\alpha \bigg(1-i\epsilon(\theta \mu(B(x,r))\tan(\pi \alpha/2)b \bigg) \\
&\hspace*{10.5cm}\times k(x,y) dx C_\beta r^{-\beta -1} dr.
\end{align*}
\\
2. From (2) in Proposition \ref{prop:comportement-general-cox}, for $\mu\in {\cal M}_{\alpha,\beta}$, we have:
\begin{align*}
&\int_{\rdr+} \psi_G\bigg(\theta \mu(B(x,r))\bigg)k_\rho(x,y) dxf\left(\frac{r}{\rho}\right) \frac{dr}{\rho}\\
&\hspace*{4.5cm}\underset{\rho\to 0}{\sim}\lambda(\rho)\rho^\beta \int_{\rdr+} \psi_G\bigg(\theta \mu(B(x,r))\bigg)k(x,y) dx C_\beta r^{-\beta -1} dr.
\end{align*}
Because $\kappa(\rho)\lambda(\rho)\rho^\beta \underset{\rho \rightarrow 0}{\longrightarrow}a \in]0,+\infty[$ and $\kappa(\rho)\longrightarrow +\infty$, necessarily we have $\lambda(\rho)\rho^\beta\longrightarrow 0$. Hence, from \eqref{eq:equiv-1-exp}, and using \eqref{hyp:f}, \eqref{hyp-krho}, we have:
\begin{align*}
&\underset{\rho\to 0}{\lim}~\left(1-\exp\left(\int_{\rdr+} \psi_G\bigg(\theta \mu(B(x,r))\bigg)k_\rho(x,y) dxf\left(\frac{r}{\rho}\right) \frac{dr}{\rho}\right) \right)\kappa(\rho)\\
&\hspace*{5.5cm}=-\int_{\rdr+}  \psi_G\bigg(\theta \mu(B(x,r))\bigg)k(x,y) dx aC_\beta r^{-\beta -1}dr.
\end{align*}
\\
3. Point (3) in Proposition \ref{prop:comportement-general-cox} writes
\begin{align*}
\int_{\rdr+} \psi_G\left(\frac{\theta \mu(B(x,r))}{n(\rho)}\right)&k(x,y) dxf\left(\frac{r}{\rho}\right) \frac{dr}{\rho}\nonumber\\
&\underset{\rho\rightarrow 0}{\sim}\frac{\lambda(\rho)\rho^\beta}{n(\rho)^\gamma}\int_{\rdr+}\psi_G(\theta  \varphi(x)r^d)k(x,y)dx v_d^{\gamma }C_\beta r^{-\beta -1}dr\nonumber\\
&=\frac{1}{\kappa(\rho)}\int_{\rdr+}\psi_G(\theta  \varphi(x)r^d)k(x,y)dx v_d^{\gamma }C_\beta r^{-\beta -1}dr
\end{align*}
if the two conditions $n(\rho)\longrightarrow 0$ and $\rho^d=o(n(\rho))$ are satisfied.\\
The first condition is clearly satisfied and for the second one we have:
\[\frac{n(\rho)^\gamma}{\rho^{\gamma d}}=\frac{\kappa(\rho)\lambda(\rho)\rho^\beta}{\rho^\beta}=\kappa(\rho)\lambda(\rho)\longrightarrow +\infty\]
which implies
\[\frac{n(\rho)}{\rho^{ d}}\longrightarrow +\infty\]
and the result follows.\\
Since $\displaystyle\underset{\rho\rightarrow 0}{\lim}~\kappa(\rho)=+\infty$, \eqref{eq:equiv-1-exp} gives, with \eqref{hyp:f} and \eqref{hyp-krho}:
\begin{align*}
&\underset{\rho\to 0}{\lim}~\left(1-\exp\left(\int_{\rdr+} \psi_G\left(\frac{\theta \mu(B(x,r))}{n(\rho)}\right)k_\rho(x,y) dxf\left(\frac{r}{\rho}\right) \frac{dr}{\rho}\right) \right)\kappa(\rho)\\
&\hspace*{6cm}=-\int_{\rdr+}\psi_G(\theta  \varphi(x)r^d)k(x,y)dx v_d^{\gamma }C_\beta r^{-\beta -1}dr. 
\end{align*}
\CQFD

\appendix

\section{Appendix: Generalities about Cox process} \label{appendix}

Let $(E,\mathcal{E})$ a measurable space. For a random measure $C$ on $(E,\mathcal{E})$, we set $\Upsilon_C$ for its characteristic function given by 
\begin{equation}\label{laplace-qcq}
\Upsilon_C(g)=\mathbb{E}\left[\exp\left(i\int_E g(x)C(dx)\right) \right]
\end{equation}
for $g:E\longrightarrow \mathbb{C}$ such that the quantity in \eqref{laplace-qcq} exists, see \cite{DVJ1}.\\
In this appendix, we specify \eqref{laplace-qcq} whose use is crucial in our argument, for Cox point process. First, we recall \eqref{laplace-qcq} for the classical case of Poisson point process on $(E,\mathcal{E})$: 

\begin{df}
Let $\lambda$ be a $\sigma$-finite measure on $(E,\mathcal{E})$ and $N$ be a point process on $E$. We say that $N$ is a Poisson point process on $E$ with intensity measure $\lambda$ if:
\begin{enumerate}
\item For all $A\in\mathcal{E}$ such that $\lambda(A)<+\infty$, $N(A)$ is a Poisson random variable with parameter $\lambda(A)$.
\item For all $n\geq 1$ and all $A_1,\dots, A_N \in\mathcal{E}$ with no intersection, $N(A_1),\dots, N(A_n)$ are mutually independent.
\end{enumerate}
\end{df}
Then \eqref{laplace-qcq} specializes as follows for a Poisson point process.
\begin{prop}\label{prop:laplace-poisson}
Let $N$ be a Poisson random measure on $E$ with intensity measure $\lambda$. Then we have:
\begin{equation*}
\ee\left[\exp\left(i\int_E g(x)N(dx)\right)  \right]=\exp \left(-\int_E \left(1-e^{i g(x)}\right)\lambda (dx) \right)
\end{equation*}
for all $g:E\longrightarrow \mathbb{C}$ such that 
\begin{equation}\label{condition-poisson}
\int_E\left(1\wedge \left|g(x)\right|\right) \lambda (dx) <+\infty.
\end{equation} 
\end{prop}

The distribution of a Poisson point process is characterized by its deterministic intensity measure as appears from Proposition \ref{prop:laplace-poisson}. A natural extension of a Poisson process is to consider a random intensity measure. From a modelling point of view, considering random intensity measure rather than deterministic ones allows to consider random constraints for the repartition of the points in the space. This forms the class of so-called Cox processes more specifically defined as follows.
\begin{df}\label{def-cox}
A point process $C$ is a Cox process directed by the random intensity function $\Lambda$ if, conditionally to $\Lambda=\lambda$, $C$ is a Poisson process with intensity measure $\lambda$.
\end{df}
In the sequel, we shall use the notation $C$ both for the locally finite collection of points $X \in C$ and for the associated random measure $\sum_{X\in C}\delta_{X}$.\\
Like for the Poisson point process, the characteristic function appears to be a suitable tool to investigate Cox process. In several particular cases, for instance in the shot-noise model, explicit expression of the characteristic function is available. \\
Now, we specify the characteristic function of a Cox process.
\begin{prop}
Let $C$ be a Cox process on $E$ directed by $\Lambda$. For all $g$ such that the quantity below is well defined we have:
\[\Upsilon_C(g)=\ee \left[ \exp\left(-\int_E \left(1-e^{ig(x)} \right)\Lambda(dx) \right) \right]. \]
\end{prop}
\Proof~ Since conditionally to $\Lambda$, $C$ is a Poisson point process with intensity $\Lambda$, Proposition~\ref{prop:laplace-poisson} ensures:
\begin{align*}
\ee \left[\exp\left.\left(i\int_E g(x)C(dx)\right) \right| \Lambda \right]=\exp\left(-\int_E \left(1-e^{ig(x)} \right)\Lambda(dx) \right)
\end{align*} 
and 
\begin{align*}
\Upsilon_C(g)=\ee \left[\exp\left(i\int_E g(x)C(dx)\right) \right]&=\ee \left[ \ee \left[\exp\left.\left(i\int_E g(x)C(dx)\right) \right| \Lambda \right]\right]\\
&=\ee \left[ \exp\left(-\int_E \left(1-e^{ig(x)} \right)\Lambda(dx) \right) \right].
\end{align*}
\CQFD

\begin{prop}\label{prop:laplace-cox-centree}
Let $C$ be a Cox process on $E$ directed by $\Lambda$. Then we have:
\begin{equation}\label{eq:laplace-centree}
\ee\left[\exp\left(i\left(\int_E g(x)C(dx)-\ee\left[\int_E g(x)C(dx)\bigg| \Lambda\right]\right)   \right) \right]=\ee \left[\exp\left(\int_E \psi(g(x))\Lambda(dx)\right) \right]
\end{equation}
where $\psi(u)=e^{iu}-1-iu$, for all $g:E\longrightarrow \mathbb{C}$ such that \eqref{eq:laplace-centree} is well defined.
\end{prop}
\Proof~ Since conditionally to $\Lambda$, $C$ is a PPP with intensity $\Lambda$ we have:
\begin{align*}
&\ee\left[\left.\exp\left(i\left(\int_E g(x)C(dx)-\ee\left[\int_E g(x)C(dx)\bigg| \Lambda\right]\right)   \right)  \right| \Lambda\right]\\
&\hspace*{5cm}=\exp\left(-i\int_E g(x)\Lambda(dx)\right)\ee\left[\left. \exp\left(i\int_E g(x)C(dx)   \right)   \right| \Lambda\right]\\
&\hspace*{5cm}=\exp\left(-i\int_E g(x)\Lambda(dx)\right)\exp\left(-\int_E \left(1-e^{ig(x)} \right)\Lambda (dx) \right)\\
&\hspace*{5cm}=\exp\left(\int_E \psi(g(x))\Lambda(dx)\right).
\end{align*}
Finally, we have:
\begin{align*}
&\ee\left[\exp\left(i\left(\int_E g(x)C(dx)-\ee\left[\int_E g(x)C(dx)\bigg| \Lambda\right]\right)   \right) \right]\\
&\hspace*{3.2cm}=\ee \left[  \ee\left[\left.\exp\left(i\left(\int_E g(x)C(dx)-\ee\left[\int_E g(x)C(dx)\bigg| \Lambda\right]\right)   \right)  \right| \Lambda\right] \right]\\
&\hspace*{3.2cm}=\ee \left[\exp\left(\int_E \psi(g(x))\Lambda(dx)\right) \right]. 
\end{align*}\CQFD

%%%%%%%%%%%%%%%%%%%%%%%%%%%%%%%%%%%%%%%%%%%%%%%%%%%%%%%%%%%%%%%%%%%%

\section*{Acknowledgement}

The author thanks his thesis advisor Jean-Christophe Breton for his various readings of preliminary versions of this paper.

%%%%%%%%%%%%%%%%%%%%%%%%%%%%%%%%%%%%%%%%%%%%%%%%%%%%%%%%%%%%%%%%%%%%

%\bibliographystyle{plain}
%\bibliography{bib.bib}

\end{document}